\documentclass[11pt]{amsart}
\usepackage{amsmath,amssymb,amscd,xypic}
\begin{document}

\renewcommand{\smallskip}{\vspace{.1in}}
\newcommand{\Sym}{{\rm Sym}}          
\newcommand{\Res}{{\rm Res}}           
\newcommand{\Hom}{{\rm Hom}}        
\renewcommand{\span}{{\rm Span}}          
\newcommand{\codim}{{\rm codim}}          
\newcommand{\RP} {\mathbb{RP}}
\newcommand{\Proj} {{\mathbb P}}
\newcommand{\N}{{\mathbb N}}
\newcommand{\R}{{\mathbb R}}
\newcommand{\A}{{\mathbf A}}
\newcommand{\bu}{{\mathbf u}}
\newcommand{\bv}{{\mathbf v}}
\newcommand{\bw}{{\mathbf w}}
\newcommand{\bx}{{\mathbf x}}
\newcommand{\bz}{{\mathbf 0}}
\newcommand{\bc}{{\mathbf c}}
\newcommand{\bfz}{{\mathbf z}}
\newcommand{\ind}{{\rm ind}}

\newtheorem{theorem}{Theorem}
\newtheorem{lemma}[theorem]{Lemma}
\newtheorem{proposition}[theorem]{Proposition}
\newtheorem{definition}[theorem]{Definition}
\newtheorem{corollary}[theorem]{Corollary}

\title{The Higher Connectivity of Intersections of Real Quadrics}
\author{Michael Larsen}
\email{larsen@math.indiana.edu}
\address{Department of Mathematics\\
    Indiana University \\
    Bloomington, IN 47405\\
    U.S.A.}

\thanks{The first named author was partially supported by NSF grant DMS-0354772.}

\author{Ayelet Lindenstrauss}
\email{ayelet@math.indiana.edu}
\address{Department of Mathematics\\
    Indiana University \\
    Bloomington, IN 47405\\
    U.S.A.}
\date{August 20,2005}
\begin{abstract}
A linear system of real quadratic forms defines a real projective variety.
The real non-singular locus of this variety (more precisely of the underlying scheme) 
has a highly connected double cover as long as each non-zero form in the system has sufficiently high Witt index.
\end{abstract}
\maketitle

There is an extensive literature (see, e.g., \cite{GM} and the references contained therein) concerning connectivity theorems for complex projective varieties.  The philosophy, ever since Lefschetz proved his hyperplane theorem, has been that complex projective subvarieties of
low codimension should have the same low dimensional homology
and homotopy groups as their ambient varieties.  The situation for
real algebraic varieties is more delicate.  Even in the simplest
non-trivial case, non-singular quadric hypersurfaces in $\RP^n$, one cannot make a connectedness statement without further hypotheses.
Indeed, pulling back the quadric
\[  Q=x_0^2+\cdots+x_{m-1}^2-x_m^2-\cdots-x_n^2=0  \]
by the double cover $S^n\to\RP^n$, we obtain a space homeomorphic
to 
\[S^{m-1}\times S^{n-m}\]
which is $r$-connected if and only if the {\em Witt index} of $Q$,
i.e. the number of mutually orthogonal hyperbolic planes in the inner
product space defined by $Q$, is at least $r+1$.  In this paper we give a qualitative generalization of this result to intersections $\Proj Y_W$
of quadric hypersurfaces.  

Such intersections were considered from a rather different point of view in \cite{IL}.  That paper proved weak 
approximation for rational points on $\Proj Y_W$ under an admissibility condition given
in Definition~\ref{admissible} below.   Admissibility guarantees
that many lines in the ambient projective space lie in $\Proj Y_W$.  In fact, for any two
points $P,Q\in\Proj Y_W$, there exists $R\in \Proj Y_W$ such that $\overline{PR}$ and
$\overline{QR}$ lie in $\Proj Y_W$.
This of
course implies that the variety is connected.  Bo-Hae Im and the first-named
author asked whether one can prove higher connectivity
results by similar methods.  In this paper we more or less carry
out that program.

To be more precise, let $V$ be an $n$-dimensional real vector space and let $W$ be a $k$-dimensional subspace of $\Sym^2 V^*$, which we view sometimes as the space of symmetric bilinear and sometimes as the space of
quadratic forms on $V$.  Let $\Proj Y_W \subset 
\Proj V\cong \RP^{n-1}$ denote the intersection of all quadric hypersurfaces given by elements of $W$.  Let $Y^{ns}_W$ denote the the set of vectors $\bv\in V$ which are not null with respect to any non-zero $\bw\in W$.  
The notation is suggested by the observation that $Y^{ns}_W$ is in fact the non-singular locus of $Y_W$ regarded as a scheme.  (This may be smaller than the non-singular locus of 
the underlying variety $Y_W$:
in the extreme case that $W$ is generated by a single quadratic form which is the square of a
non-zero linear form on $V$, the scheme has no smooth points at all, but the underlying variety is just a hyperplane.)  Obviously $Y^{ns}_W$ is closed under scalar multiplication by positive real numbers, so it makes sense to talk about $\Proj Y^{ns}_W$.  

\begin{definition}
\label{admissible}
We say that $W$ is {$m$-admissible} for some $m\in \N$ 
if for every $\bw\in W$, $\bw$ is positive definite on
a subspace of $V$ of dimension greater than or equal to $m$ and negative definite on a subspace of $V$ of dimension greather than or equal to $m$.
\end{definition}

We can now state our basic result.

\begin{theorem} 
There exists a function $m\colon\N^2\to\N$ such that if $i\in \N$, $0<k\in\N$
and $W$ is an $m(i,k)$-admissible $k$-dimensional space of 
quadratic forms on $V$ then $\Proj Y^{ns}_W$ has a unique
$i$-connected double cover.  
\label{IntroThm}
\end{theorem}

We remark that generically, $\Proj Y_W$ is a non-singular projective variety, and in this case,
$\Proj Y^{ns}_W = \Proj Y_W$.

Theorem~\ref{IntroThm} follows directly from Theorem \ref{MainTheorem} below; the latter is slightly more general and also gives an explicit formula for the function $m$.
The idea of the proof is to contract an arbitrary continuous map
$f\colon S^i\to \Proj Y^{ns}_W$  to a point $P$ along a cone connecting $P$
to $f(S^i)$.  Unfortunately, we cannot expect that there exists
a point $P\in \Proj Y^{ns}_W$ for which all the lines from $p$ to
$f(S^i)$ lie in $\Proj Y^{ns}_W$.  We might try to achieve the contraction in
two steps by fixing $P$ and constructing
$g\colon S^i\to \Proj Y^{ns}_W$ such that for every $x\in S^i$, there
are lines in $\Proj Y^{ns}_W$ connecting $P$ with $g(x)$ 
and $g(x)$ with $f(x)$.  This can be regarded as a lifting problem
\[\xymatrix{&Z\ar[d] \\ 
    S^i\ar@{-->}[ur]\ar[r]^f& \Proj Y_W^{ns}     
}\]
where $Z$ consists of ordered pairs $(Q,R)\in (\Proj Y^{ns}_W)^2$
for which the lines $\overline{PQ}$ and 
$\overline{QR}$ lie in $\Proj Y^{ns}_W$.  The fibers of $(Q,R)\mapsto R$ are essentially intersections of quadrics of lower
admissibility than the original system.

If $Z\to \Proj Y^{ns}_W$
were a fibration, one could use the long exact homotopy sequence
to deduce the existence of a lifting from the $(i-1)$-connectivity of the fibers.  Unfortunately, this is not the case.
Instead, we make use of a similar lifting theorem for submersions of smooth manifolds where
the inverse image of every point is $(i-1)$-connected.  All the complications in our argument are caused by our need to pass from 
the singular spaces $\Proj Y^{ns}_W$ and $Z$ to open subsets on
which this lifting theorem can be applied.

We regard our theorem as a first illustration, in a special but nevertheless non-trivial setting, 
of the idea that real
varieties satisfying sufficiently strong rational connectedness
properties in the sense of \cite{KMM} must in fact be highly connected in the homotopy theory sense as well.  As far as we know, a general definition of higher rational connectedness has yet to 
be given, although there has been recent work in this direction \cite{HS}.

\smallskip
Throughout this paper, we identify real varieties with their real loci.  This should not cause any confusion.
\smallskip

We thank the Hebrew University of Jerusalem for its hospitality while this work was being 
carried out.

\begin{definition}
Given a vector space $V$, a subspace $V_0\subset V$, and a space of quadratic forms $W\subset \Sym^2 V^*$, the \emph{restriction of $W$ to $V_0$}, denoted $\Res_{V_0} W$,
is the image of $W$ under the 
natural surjective linear transformation 
$\Sym^2 V^* \to \Sym^2V_0^*$ to $W$. 
\end{definition}

We start with the following lemma (\cite[Lemma 3]{IL}),
\begin{lemma}
Let $V$ be a finite-dimensional real vector space, and let $V_0\subset V$ be a subspace of codimension $d$.  Let 
$W$ be an $m$-admissible, $k$-dimensional space of 
quadratic forms on $V$.  Then if $W$ is $m$-admissible for some
$m>d$, $Res_{V_0} W$ is $(m-d)$-admissible.
\label{ReducingLemma}
\end{lemma}

\begin{proof}
Since $W$ is $m$-admissible on $V$, for each $\bw\in W$ there
are subspaces $V^+_\bw$ and $V^-_\bw$ of $V$ of dimension $m$ 
each on which $\bw$ is positive definite and negative definite,
respectively.  Then $V_0\cap V^+_\bw$ and $V_0\cap V^-_\bw$
both have dimension greater than or equal to $m-d$, and $\bw$ is
positive definite and negative definite, respectively, on them.
\end{proof}

According to \cite[Prop. 4]{IL}, if the $k$-dimensional space of quadratic forms on 
$V$ is $(k^2-k+1)$-admissible, the (nonlinear)
evaluation map $E\colon V\to W^*$ given by 
\begin{equation} E(\bv)(\bw)=\bw(\bv,\bv)
\label{Edef}\end{equation}
is surjective.  
We show here that  we can guarantee surjectivity even after omitting a finite number of subvarieties of $V$ if their codimension is high enough.

\begin{proposition}
Let $V$ be a finite dimensional real vector space, and let
$W\subset \Sym^2 V^*$ be a $k$-dimensional, $(k^2+k-1)$-admissible  space of quadratic
forms on $V$.  Let $X\subset V$ be a subvariety, or a finite union
of subvarieties, of codimension $\geq k^2+4k-3$.  Then the evaluation
map  $E\colon V\setminus X\to W^*$ as defined in  (\ref{Edef})
above is surjective.
\label{EvalSurjProp}
\end{proposition}

\begin{proof}
At the cost of decreasing the codimension by $1$, we may assume that $X$ is a union of lines through the origin.  As $E(a\bv)=a^2E(\bv)$,
$E(V\setminus X)$ is then closed under multiplication by strictly positive
real numbers.  Also, $E(V\setminus X)$ is not contained in any closed
half-space through the origin in $W^*$ since that would imply the existence of a $\bw\in W$ with $\bw(\bv,\bv)\geq 0$ for {\em all}
$\bv\in V\setminus X$, whereas by our assumption, $\bw$ is negative
definite on some subspace so $\bw(-,-)$ is strictly negative
on a non-empty open subset of $V$ (and $X$ has an empty interior).

\smallskip

If $k=1$, $E(V\setminus X)$ is by what we have proved closed under multiplication by positive scalars, and not contained in a half-line.
So all that remains to be proved is that $\bz\in E(V\setminus X)$.
But since $W=\R\cdot \bw$ for a single $\bw\in W$, 
$\{\bv\in V\mid E(\bv)=\bz\}$ is a conic with manifold points, i.e. points whose neighborhoods are manifolds of codimension $1$, so
the conic cannot be contained entirely in $X$, which is of codimension $\geq 1^2+4\cdot 1-3=2$, and there is
$\bv\in V\setminus X$ with $E(\bv)=\bz$.

\smallskip
If $k>1$, we use induction on $k$ and 11.7.3 in \cite{Ro}, which says
that a convex subset of a finite dimensional vector space which is
not contained in any closed half-space through the origin, and is
closed under multiplication by positive scalars must be the entire
vector space.  
To use this result, we still need to prove that $E(V\setminus X)$ is
convex.  Consider $\bv_1, \bv_2\in V\setminus X$.  We need to prove
that the line segment between $E(\bv_1)$ and $E(\bv_2)$ is contained in $E(V\setminus X)$.  

If $E(\bv_1)$ and $E(\bv_2)$ lie on the same closed ray through the
origin, this follows trivially from the closure of $E(V\setminus X)$
under multiplication by positive scalars.

Otherwise, we will need to use the induction hypothesis.  Set 
\[W_1=\bigl( W^*/\R\cdot E(\bv_1)\bigr)^*\subset W.\]  
Since we know that $E(\bv_1)\neq\bz$ (otherwise we would be in the first case), $W_1$ is $(k-1)$-dimensional.  Now let
\[V_{1,2}=\{\bv\in V\mid \bv\perp_W\bv_1\  {\rm and}\  \bv\perp_W\bv_2\}.  \]
where we use the notation $\bv\perp_W\bu$ to indicate that $\bw(\bv,\bu)=0$ for all $\bw\in W$.
Then $V_{1,2}$ has codimension $\leq 2k$ in $V$.  The restriction map $W_1\to\Res_{V_{1,2}}W_1$ is
an isomorphism, since for all $\bw\in W_1\subset W$, if $\bw$ restricts to the zero form
in $\Sym^2 V_{1,2}^*$, it must have a nullspace of codimension $\leq 2k$, whereas we know
that it is positive definite and negative definite on subspaces of dimension $\geq k^2+k-1$ each.

Since the restriction of $W_1$ to $V_{1,2}$ is an isomorphism and $\codim V_{1,2}\leq 2k$, by Lemma
\ref{ReducingLemma},  $\Res_{V_{1,2}}W_1$ is $(k^2-k-1)$-admissible.  We can write
$k^2-k-1=(k-1)^2+(k-1)-1$.  

Let $X_{1,2}$ denote the Zariski-closure of  
\[\span(\bv_1,\bv_2,X)\cap V_{1,2}=\bigcup_{\bx\in X}\span(\bv_1, \bv_2, \bx)\cap V_{1,2}.\]
As $\span(\bv_1,\bv_2,X)$ is the image of a morphism $\A^3\times X\to V$, it follows that
\begin{multline*}
\codim_{V_{1,2}} X_{1,2} \geq\codim_V X -3-2k\\  \geq
k^2+4k-3-3-2k=(k-1)^2+4(k-1)-3.\\
\end{multline*}

By the induction hypothesis, $E(V_{1,2}\setminus X_{1,2})=W_1^*$, so we can find 
$\bv_3\in V_{1,2}\setminus X_{1,2}$ with 
$E(\bv_3)=E(\bv_2)$ in $W^*_1\cong W^*/\R\cdot E(\bv_1)$.
This means that
\begin{equation}
E(\bv_3)=E(\bv_2)+c E(\bv_1)
\label{vthreedef}
\end{equation}
in $W^*$.  Since $\bv_3\in V_{1,2}\setminus X_{1,2}$, we know the following:

\begin{enumerate}
\item  $\span(\bv_1,\bv_3)\cap X =\{0\}$ (otherwise $\bv_3\in\span(\bv_1, X)$);
\item $\span(\bv_2,\bv_3)\cap X =\{0\}$ (otherwise $\bv_3\in\span(\bv_2, X)$);
\item $\bv_3\perp_W\bv_1$, so $E(a\bv_1+b\bv_3)=a^2 E(\bv_1)+b^2 E(\bv_3)$ for all $a,b\in \R$;
\item  $\bv_3\perp_W\bv_2$, so $E(a\bv_2+b\bv_3)=a^2 E(\bv_2)+b^2 E(\bv_3)$ for all $a,b\in \R$.
\end{enumerate}

First assume that $E(\bv_1)$ and $E(\bv_2)$ lie on the same line through the origin.  We have already
dealt with the easy case in which they lie on the same closed ray through the origin, and explained why,
if they sit on opposite rays, the open rays they sit on are in $E(V\setminus X)$.  The only issue, then,
is whether the origin is in $E(V\setminus X)$.  If $E(\bv_3)=\bz\in W^*$, we are done. 
If
$E(\bv_3)$ sits on the same open ray as $E(\bv_2)$, look at the path 
$\{E(t\bv_1+(1-t)\bv_3)\mid t\in[0,1]\}$
from $E(\bv_3)$ to $E(\bv_1)$.  This path goes through the origin, because $E(\bv_3)$ and $E(\bv_1)$ lie
on opposite rays, and by Property 3 above, the image of the path is a line segment.  
The path lies in $E(V\setminus X)$ because it is the image under $E$ of the
path  $\{t\bv_1+(1-t)\bv_3\mid t\in[0,1]\}$, which by Property 1 above  could a priori intersect $X$ only at $\bz\in V$,
and it is impossible for $\bz$ to be written as  $t\bv_1+(1-t)\bv_3$ for any $t$ since that would imply that
$\bv_1$ and $\bv_3$ are linearly dependent, which would make $E(\bv_1)$ and $E(\bv_3)$ sit
on a single closed ray through the origin, rather than on opposite open rays.  We get $\bz\in E(V\setminus X)$.
Finally, if $E(\bv_3)$ sits on the same open ray as $E(\bv_1)$, repeat the last argument with 
$\bv_2$ in the role of $\bv_1$, to obtain as before $\bz\in E(V\setminus X)$.

\smallskip

So assume now that $E(\bv_1)$ and $E(\bv_2)$ do not sit on a common line through the origin.
Observe from equation (\ref{vthreedef}) that $E(\bv_1)$, $E(\bv_2)$, and $E(\bv_3)$ all lie on
a single plane through the origin, and $E(\bv_2)$ and $E(\bv_3)$ lie on the same side of the
line through $E(\bv_1)$ on that plane.    

By the same equation (\ref{vthreedef}),  neither $E(\bv_1)$ and $E(\bv_3)$
nor $E(\bv_2)$ and $E(\bv_3)$ sit on a single closed ray through the origin.  So, as  we argued before,
the paths  $\{E((1-t) \bv_1+t\bv_3)\mid t\in[0,1]\}$ and $\{E((1-t) \bv_3+t\bv_2)\mid t\in[0,1]\}$
connecting $E(\bv_1)$ to $E(\bv_3)$ and then to $E(\bv_2)$ by line segments are
all in $E(V\setminus X)$.  

Now we look at the plane on which the three vectors lie, all on a single half-plane, and observe 
that the line segment from $E(\bv_1)$ to $E(\bv_2)$ lies between the origin and the line segments
from $E(\bv_1)$ to $E(\bv_3)$ and from  $E(\bv_3)$ to $E(\bv_2)$, that is: each point on the 
line segment from $E(\bv_1)$ to $E(\bv_2)$ is a multiple by a strictly positive real number of
a point in the broken path  from  $E(\bv_1)$ to $E(\bv_3)$ and then to $E(\bv_2)$.  Therefore 
the line segment from $E(\bv_1)$ to $E(\bv_2)$ is contained in $E(V\setminus X)$.
\end{proof}

Proposition \ref{EvalSurjProp} above will be the key technical ingredient in proving that certain
maps are submersions.

\begin{definition}
Let $V$ be a vector space with $W\subset \Sym^2 V^*$ a $k$-dimensional space of quadratic forms
on $V$.  We can map $\phi\colon V\otimes W\to V^*$ by setting $\phi(\bv\otimes \bw)= \bw(\bv, -)$.  We say that
an $n$-tuple $(\bv_1,\ldots,\bv_n)\in V^n$ is \emph{$W$-independent} if $\phi(\span (\bv_1,\ldots,\bv_n)
\otimes W)\subseteq V^*$ is $nk$-dimensional.
\end{definition}

Observe that the set of $W$-independent $n$-tuples is open in $V^n$.  Indeed, if we fix a positive definite inner product on $W$ and let $S(W^n)$ denote the (compact) set of unit vectors, then the complement of the set of $W$-independent $n$-tuples is the (proper) projection to $V^n$ of \[\{((\bv_1,\ldots,\bv_n),(\bw_1,\ldots,\bw_n))\in V^n\times S(W^n)\mid 
\phi(\bv_1\otimes\bw_1+\cdots+\bv_n\otimes\bw_n) = \bz\}\]

We now introduce notation which will be used for the rest of the paper.
Let $V^n_{\ind}$ be the set of all $(\bv_1,\ldots,\bv_n)$ in $V^n$  which are $W$-independent.
Let $X_n$ consist of all $(\bv_0, \bv_1,\ldots,\bv_n)\in V^{n+1}_{\ind}$  such that $\bv_0\perp_W\bv_i$
for all $0\leq i\leq n$.  Projection on the last $n$ coordinates gives a natural map
\begin{equation}
\label{XY}
p\colon X_n\to V^n_{\ind}
\end{equation}

\begin{lemma}
The variety $X_n$ is a manifold. 
\end{lemma}

\begin{proof}
We will look at the function $\theta\colon V^{n+1}\to (W^*)^{n+1}$ defined by
\begin{equation}
\label{theta}
\theta(\bv_0,\bv_1,\ldots,\bv_n)=
\{\bw\mapsto(\bw(\bv_0,\bv_0),\bw(\bv_0,\bv_1),\ldots, \bw(\bv_0,\bv_n))\} .
\end{equation}
Then $X_n$ is the intersection of $\theta^{-1}(\bz)$ with the open set $V^{n+1}_{\ind}$.  But on the
set  $V^{n+1}_{\ind}$  of $W$-independent $(n+1)$-tuples, 
$\theta_*\colon T(V^{n+1})\to T((W^*)^{n+1})$
is surjective.  To see this, recall from the definition of $W$-independence that it requires
$\{\bv_0,\bv_1,\ldots,\bv_n\}$ to be linearly independent and the map 
\[\phi\colon \span(\bv_0,\bv_1,\ldots,\bv_n)\otimes W\to V^*\]
given by $\phi(\bv\otimes\bw)=\bw(\bv,-)$ to be injective.  This implies that the dual map
\begin{multline*}
V\cong V^{**}\xrightarrow{\phi^*}
\span(\bv_0,\bv_1,\ldots,\bv_n)^*\otimes W^*\cong \Hom(\span(\bv_0,\bv_1,\ldots,\bv_n), W^*)\\
\cong\bigoplus_{i=0}^n\Hom(\R\cdot \bv_i, W^*)\cong (W^*)^{n+1}\\
\end{multline*}
given by $\phi^*(\bv)=\{\bw\mapsto(\bw(\bv,\bv_0),\bw(\bv,\bv_1), \ldots,\bw(\bv,\bv_n))\}$ is surjective.

Now at any $\bv  := (\bv_0,\bv_1,\ldots,\bv_n)\in V^{n+1}$ and
any $(\bu_0,\bu_1,\ldots,\bu_n)\in T_{\bv}(V^{n+1})$,
\begin{multline*}
\theta_*(\bu_0,\bu_1,\ldots,\bu_n)\\
=\lim_{\epsilon\to 0}\frac{\theta(\bv_0+\epsilon\bu_0,\bv_1+\epsilon\bu_1,\ldots, \bv_n+\epsilon\bu_n)-\theta(\bv_0,\bv_1,\ldots\bv_n)}{\epsilon} \\
=\{\bw\mapsto(2\bw(\bu_0,\bv_0), \bw(\bv_0,\bu_1)+\bw(\bu_0,\bv_1),\ldots,\bw(\bv_0,\bu_n)+\bw(\bu_0+\bv_n))\}.
\end{multline*}
We can look only on those tangent vectors $(\bu_0,\bu_1,\ldots,\bu_n)$ where  $\bu_1=\bu_2=\cdots
=\bu_n=\bz$ and still get surjectivity:  for these,
\[\theta_*(\bu_0,\bz,\ldots,\bz)=\{\bw\mapsto (2\bw(\bu_0,\bv_0), \bw(\bu_0,\bv_1),\ldots,\bw(\bu_0,\bv_n))\},\]
which differs from $\phi^*$, known to be surjective on $V^{n+1}_{\ind}$, by the invertible transformation 
of scaling the first coordinate by a factor of $2$.

It remains only to invoke the implicit function theorem for $\theta\vert_{V_{\ind}^{n+1}}$, and we get that
$X_n=\theta\vert_{V_{\ind}^{n+1}}^{-1}(\bz)$ is a manifold.
\end{proof}

\begin{corollary}
The variety $Y_W^{ns}$ is an open manifold of codimension $k$ (the dimension of $W$) in $V$.
\label{DimYW}
\end{corollary}

\begin{proof}
By our definitions, $Y_W^{ns}=X_0$ is the fiber, by the previous lemma, of a submersion of manifolds
$V_{\ind}\to W^*$ so its dimension is $\dim V-\dim W^*=\dim V-k$.
\end{proof}

\begin{proposition}
Let $W\subset\Sym^2 V^*$ be a $k$-dimensional space of quadratic forms on $V$ which
is $M$-admissible, where $M=\lceil\frac{k^2+3nk+5k-3}{2}\rceil$.  Then the map
$p\colon X_n\to V^n_{\ind}$ of (\ref{XY}) is surjective and submersive.
\label{FirstLemma}
\end{proposition}

\begin{proof}
We start by showing that $p$ is surjective.  Given $(\bv_1,\ldots,\bv_n)\in V^n_{\ind}$, we need to find
$\bv_0$ such that $\bv_0\perp_W\bv_0$ and $\bv_0\perp_W\bv_i$, $1\leq i \leq n$, and 
$\{\bv_0,\bv_1,\ldots,\bv_n\}$ are $W$-independent.  Since we know that $\{\bv_1,\ldots,\bv_n\}$ 
are $W$-independent, this last condition translates to having $\bv_0$ not be a null vector of
any $\bw\in W$ and having 
\[\phi(\bv_0\otimes W)\cap \phi(\span(\bv_1,\ldots,\bv_n)\otimes W)=\bz.\]
Let $V'=\{\bv\in V \mid \bv\perp_W\bv_i,\ 1\leq i \leq n\}$. By the $W$-independence, $V'$ is of codimension $nk$ in $V$.  Since $W$ is $M$-admissible on $V$, by Lemma \ref{ReducingLemma}, the restriction
of $W$ to $V'$, which we will also denote $W$, is $(M-nk)$-admissible.

What we need to do is to find $\bv_0\in V'$ with $E(\bv_0)=\bz$  which 
does not belong to
\[X'=\{\bv_0\in V'\mid \exists \bz\neq\bw\in W,\ \bc\in\span(\bv_1,\ldots,\bv_n)\otimes W:\ \phi(\bv_0\otimes\bw)=\bc\}.\]
We want to show that the codimension of $X'$  is high enough that we can use Proposition
\ref{EvalSurjProp}.   We consider the quasi-affine subvariety 
\[Z' = \{(\bv_0,\bw,\bc)\!\in V'\times W\times\phi(\span(\bv_1,\ldots,\bv_n)\otimes W)
\mid \bw\neq\bz,\ \phi(\bv_0\otimes\bw)=\bc\}.\]
Projecting to the second and third coordinates
\[Z'\to W\times\phi(\span(\bv_1,\ldots,
\bv_n)\otimes W)\setminus \{\bz\}\times\phi(\span(\bv_1,\ldots,
\bv_n)\otimes W),\]
the fiber over each point $(\bw,\bc)$ is the set of all $\bv_0$ such that $\phi(\bv_0\otimes\bw)=\bc$,
which, if non-empty, has the same dimension as the set of all $\bv$ such that $\phi(\bv\otimes\bw)=\bz$.  By our admissibility condition, this set has codimension $\geq2(M-nk)$ (since vectors in
the positive definite subspace or in the negative definite subspace cannot be null vectors).  So
\begin{multline*}
\dim Z'\leq\dim W+\dim(\phi(\span(\bv_1,\ldots,\bv_n)\otimes W))
+(\dim V'-2(M-nk))
\\=k+ nk+\dim V'-(2(M-nk))=\dim V'-(2M-3nk-k)
\end{multline*}
(see, e.g., \cite[II~Ex.~3.22(b)]{Har}.)
Now $X'$ is the image of $Z'$ under  projection to the first coordinate so 
\[\dim X'\leq\dim Z'\]
and
\begin{multline*}
\codim X'\ \geq\  2M-3nk-k\
\geq k^2+3nk+5k-3-3nk-k=k^2+4k-3
\end{multline*}
as required by Proposition \ref{EvalSurjProp}.   Since the conditions of Proposition \ref{EvalSurjProp} hold, there exists $\bv_0\in V'\setminus X'$
with $E(\bv_0)=\bz$, and then $(\bv_0,\bv_1,\ldots,\bv_n)\in X_n$ and $p(\bv_0,\bv_1,\ldots,\bv_n)=
(\bv_1,\ldots,\bv_n)$.

\smallskip

Now we need to show that $p$ is a submersion.  Given  $(\bv_0,\bv_1,\ldots,\bv_n)\in X_n$,
we have to show that for any vector $(\bu_1,\ldots,\bu_n)\in V^n\cong T(V^n_{\ind})$, there is some
$\bu_0\in V$ so that $(\bu_0,\bu_1,\ldots, \bu_n)$ is tangent to $X_n$ at $(\bv_0,\bv_1,\ldots,\bv_n)$.

Recall that for the function $\theta\colon V^{n+1}\to (W^*)^{n+1}$ given in equation (\ref{theta}),
$X_n=V^{n+1}_{\ind}\cap \theta^{-1}(\bz)$.  So $(\bu_0,\bu_1,\ldots,\bu_n)\in V^{n+1}$
is tangent to $X_n$ at  $(\bv_0,\bv_1,\ldots,\bv_n)\in V^{n+1}_{\ind}$ if and only if 
$\theta(
\bv_0+\epsilon\bu_0,
\bv_1+\epsilon\bu_1,\ldots
\bv_n+\epsilon\bu_n)$
is $O(\epsilon^2)$ as $\epsilon\to 0$.  Now
\begin{multline*}
\theta(
\bv_0+\epsilon\bu_0,
\bv_1+\epsilon\bu_1,\ldots
\bv_n+\epsilon\bu_n)   
\\=
\{\bw\mapsto \theta(\bv_0,\bv_1,\ldots,\bv_n)(\bw)
\\+
 \epsilon
(2\bw(\bu_0,\bv_0),
\bw(\bv_0,\bu_1)+\bw(\bu_0+\bv_1),
\ldots,
\bw(\bv_0,\bu_n)+\bw(\bu_0+\bv_n))\\+
\epsilon^2  \theta(\bu_0,\bu_1,\ldots,\bu_n)(\bw)\}\\
\end{multline*}
and since $(\bv_0,\bv_1,\ldots,\bv_n)\in X_n$, $ \theta(\bv_0,\bv_1,\ldots,\bv_n)(\bw)=\bz$.  So the 
tangency condition translates to the equations
\begin{equation}
\begin{cases}
\bw(\bu_0,\bv_0)=0 & \forall\bw\in W\\
\bw(\bv_0,\bu_i)+\bw(\bu_0,\bv_i)=0 &\forall\bw\in W,\ 1\leq i\leq n.\\
\end{cases}
\label{system}
\end{equation}
Recall that $(\bv_0,\bv_1,\ldots,\bv_n)$ and $(\bu_1,\ldots,\bu_n)$ are all given, and the question
is whether we can find $\bu_0$ to satisfy the system in (\ref{system}).  Pick a basis $\bw_1,\ldots,\bw_k$
of $W$.  Then (\ref{system}) is equivalent to the inhomogeneous system of $(n+1)k$ linear 
equations on $\bu_0$
\begin{equation}
\begin{cases}
\bw_j(\bu_0,\bv_0)=0 & 1\leq j\leq k\\
\bw_j(\bu_0,\bv_i)= -\bw_j(\bv_0,\bu_i) &1\leq j\leq k,\ 1\leq i\leq n.\\
\end{cases}
\label{systemII}
\end{equation}
However the corresponding homogeneous system
\[ \bw_j(-,\bv_i)=0,\ \ \ 1\leq j\leq k,\ 0\leq i\leq n\]
is the same as requiring that $\phi(\span(\bv_0,\bv_1,\ldots,\bv_n)\otimes W)\subset V^*$ should
vanish on the vector that we choose, and the $W$-independence of $\{\bv_0,\bv_1,\ldots,\bv_n\}$
says exactly that the homogeneous system (\ref{systemII}) is $(n+1)k$-dimensional.
So any corresponding inhomogeneous system like (\ref{system}) is solvable, and we can
find a $\bu_0$ as desired.
\end{proof}

\begin{definition}
Let $Z_n$ consist of all $(\bv_0, \bv_1,\ldots,\bv_n)\in V^{n+1}$  such that $\bv_0\in Y_W^{ns}$ and
\[\phi(\{\bv_0\}\otimes W)\cap \phi(\span(\bv_1,\ldots,\bv_n)\otimes W)=\{\bz\}.\]
\end{definition}
In other words, we do not insist on 
$\{ \bv_1,\ldots,\bv_n\}$ being $W$-independent, but we do not want adding
 $\bv_0$ to add to the $W$-dependence.

\begin{lemma}
The variety $Z_n$ is a manifold, and the projection to the last $n$ coordinates
$Z_n\to V^n$ is a surjective submersion.
\label{SecondLemma}
\end{lemma}

\begin{proof}
We begin as before with verifying the surjectivity.  Consider the map 
$\gamma\colon V^{n+1}\to W^*$ given by 
\[\gamma(\bv_0,\bv_1,\ldots,\bv_n)=\{\bw\mapsto\bw(\bv_0,\bv_0)\}.\]
Then $Z_n$ is the intersection of  $\gamma^{-1}(\bz)$, $V^1_{\ind}\times V^n$, and
\[\{(\bv_0,\ldots,\bv_n)\mid
\phi(\{\bv_0\}\otimes W)\cap \phi(\span(\bv_1,\ldots,\bv_n)\otimes W)=\{\bz\}\}.\]
The two last sets are open in $V^{n+1}$,
so $Z_n$ is open in $\gamma^{-1}(\bz)$, which is a manifold because 
$\gamma_*\colon T(Z_n)\to T(W^*)$
is surjective.  This is true because for any $\bv = (\bv_0,\bv_1,\ldots,\bv_n)\in Z_n$ and 
$(\bu_0,\bu_1,\ldots,\bu_n)\in T_{\bv}(Z_n)$,
\begin{multline*}
\gamma_*(\bu_0,\bu_1,\ldots,\bu_n) \\
=\lim_{\epsilon\to 0}
\frac{\gamma(
\bv_0+\epsilon\bu_0,
\bv_1+\epsilon\bu_1,\ldots,
\bv_n+\epsilon\bu_n)-\gamma(\bv_0,\bv_1,\ldots,\bv_n)}{\epsilon}\\
=\{\bw\mapsto2\bw(\bu_0,\bv_0)\},
\end{multline*}
which is $k$-dimensional as $\bu_0$ ranges over $V$ exactly because $\bv_0$ is not a null
vector of any $\bw\in W$.  

To see that the projection to the last $n$ coordinates $Z_n\to V^n$ gives a surjection on tangent spaces, observe that $(\bu_0,\bu_1,\ldots,\bu_n)\in T_{\bv}(Z_n)$
exactly if $\bu_0\perp_W\bv_0$, which poses no restrictions on $\{\bu_1,\ldots,\bu_n\}$.
\end{proof}

Now we have all the ingredients we need to prove our main theorem, which implies Theorem \ref{IntroThm} in the introduction (after dividing by the action of $\R^\times$):
\begin{theorem}
Let $i \geq -1$ and $k> 0$ be integers, and let
\begin{multline*}
r(i,k)=k^2+2ik+i+6k-2,\ \ \ 
m(i,k)=k^2+2ik+3k+2.
\end{multline*}
Let $V$ be a finite-dimensional real vector space and $W\subset\Sym^2 V^*$
a $k$-dimensional $m(i,k)$-admissible space of quadratic forms.
Let $X$ denote a finite union of closed subvarieties of $V$ of codimension 
$\geq r(i,k)$, and suppose that $X$ contains all the null vectors of $V$, i.e., all
vectors in $V\setminus V^1_{\ind}$.
Then $Z:= Y_W\setminus X$ is $i$-connected.
\label{MainTheorem}
\end{theorem}
Note that the set of all null vectors 
\[\{\bv\in V\mid \phi(\bv\otimes\bw)=\bz \ {\rm for\ some}\ \bz\neq\bw\in W\}\]
 is exactly the union of all the nullspaces (which have codimension
$\geq 2m(i,k)$ each) of all the nonzero $\bw\in W$,
so it has codimension $\geq 2m(i,k)-k$.  But 
\[2m(i,k)-k=2k^2+4ik+5k+4>k^2+2ik+i+6k-2=r(i,k)\]
for all $i\geq -1$, $k> 0$ because for $i=-1$, $k^2-3k+7>0$ and for $i\geq 0$,
$i(2k-1)\geq 0$ and $k^2-k+6>0$ for all $k>0$.
So adding the set of null-vectors of all the $\bz\neq\bw\in W$ to $X$, if they were not already inside 
$X$, will not push the codimension of $X$ below $r(i,k)$.  

\begin{proof}
We will use induction on $i$.  In the base case $i=-1$, what we need to show is that $Z$ is non-empty.  To get $Z$, we remove from $Y_W$ the set $X$, which is a finite union of closed subvarieties
of codimension $\geq r(-1,k)=k^2+4k-3$.  By Proposition \ref{EvalSurjProp}, if $W$ is $(k^2+k-1)$-admissible, so certainly if it is $m(-1,k)=(k^2+k+2)$-admissible, $Y_W\setminus X\neq \emptyset$.

\smallskip

Now to do the inductive step, pick a basepoint $\bfz_0\in Z$, and define the following sets:
\[A_{0,4}=\{(\bv_0,\bv_4)\in Z^2\}\]
\[A_{0,2,4}=\{(\bv_0,\bv_2,\bv_4)\in Z^3\mid \phi(\bv_2\otimes W)\cap \phi(\span(\bv_0,\bv_4)\otimes W)=\{\bz\}\}\]
\begin{multline*}
A_{0,1,2,4}=\{(\bv_0,\bv_1,\bv_2,\bv_4)\mid (\bv_0,\bv_2,\bv_4)\in A_{0,2,4},\ \bv_1\perp_W\span(\bv_0,\bv_2),
\\
\bv_1\in Y_W,\ {\rm the\ lines\ from\ }\bv_1{\rm\ to\ }\bv_0{\rm\ and\ }\bv_2\ {\rm miss}\ X\}
\end{multline*}
\begin{multline*}
A_{0,1,2,3,4}=\{(\bv_0,\bv_1,\bv_2,\bv_3, \bv_4)\mid (\bv_0,\bv_1,\bv_2,\bv_4)\in A_{0,1,2,4},\ 
\bv_3\in Y_W,\\ \bv_3\perp_W\span(\bv_2,\bv_4),
\ {\rm the\ lines\ from\ }\bv_3{\rm\ to\ }\bv_2{\rm\ and\ }\bv_4\ {\rm miss}\ X\}.
\end{multline*}
We have obvious projections
\[A_{0,1,2,3,4}\to A_{0,1,2,4}\to A_{0,2,4}\to A_{0,4}.\]
The strategy is the following: given a map $f\colon S^i\to Z$, we need to show that $f$ is null-homotopic
in $Z$.  We define a map $(f_0,f_4)\colon S^i\to A_{0,4}$ where $f_0$ is the constant map at $\bfz_0$ and $f_4=f$.
Then we lift to a map $S^i\to A_{0,2,4}$ whose projection to $A_{0,4}$ is homotopic to
$(f_0,f_4)$, and continue lifting (up to homotopy) all the way up to $A_{0,1,2,3,4}$, that is:
we get a map
\[(g_0,g_1,g_2,g_3,g_4)\colon S^i\to A_{0,1,2,3,4}\]
with $f_0\simeq g_0$ and $g_4\simeq f_4=f$ as maps $S^i\to Z$.
Once we get this lifting, we are done, because by construction $g_0$ is homotopic to $g_1$ which
is homotopic to $g_2$ in $Z$, and $g_2$ is homotopic to $g_3$ which is homotopic to $g_4$ in
$Z$.  To go from $g_0$ to $g_1$, for example, we look at the homotopy 
$H\colon S^i\times[0,1]\to Z$
given by
\[H(s,t)=(1-t)g_0(s)+tg_1(s).\]
This is clearly a continuous map into $V$, but in fact it lands in $Z$: recall that $g_0(s),g_1(s)\in Z\subseteq Y_W$ for all $s$, and by construction $g_1(s)\perp_W g_0(s)$ for all $s$.
This means that for any $\bw\in W$,
\[\bw(H(s,t),H(s,t))=(1-t)^2\bw(g_0(s), g_0(s))+t^2\bw(g_1(s), g_1(s))=0\]
for all $s,t$, i.e. $H(s,t)\in Y_W$ for all $(s,t)\in S^i\times [0,1]$.  But also, since this homotopy
is along straight lines and happens in $A_{0,1,2,3,4}$, the assumption about the lines from
$\bv_1$ to $\bv_0$ and $\bv_2$ missing $X$ tells us that  $H(s,t)\in Z$ for all $(s,t)\in S^i\times [0,1]$. 

We repeat the same construction and argument for homotopies from $g_1$ to $g_2$, from $g_2$ to
$g_3$, and from $g_3$ to $g_4$.  We get that for the constant map $f_0$,
\[f_0\simeq g_0\simeq g_1\simeq g_2\simeq g_3\simeq g_4\simeq f_4=f.\]
The idea of the homotopy $g_0\simeq g_4$ is to find a function $S^i\to Z$ which is pointwise
$W$-orthogonal both to $g_0$ and to $g_4$, as well as to itself.  But, as will be explained below, 
to do that we need $g_0$ and $g_4$ to be pointwise $W$-independent, which is not always
the case: a priori $g_4$ could, for example,  be a space-filling curve which passes through every possible
candidate for $\bfz_0$. So we pick $g_2$ which {\em is} pointwise $W$-independent of
$g_0$ and also of $g_4$, and then find an intermediate $g_1$ between $g_0$ and $g_2$,
which is pointwise $W$-orthogonal to those two, and an intermediate $g_3$ between
$g_2$ and $g_4$, pointwise $W$-orthogonal to them.

\smallskip

To lift $(f_0,f_4)$, up to homotopy, from $A_{0,4}$ to $A_{0,2,4}$ we let 
\begin{multline*}
B_{0,2,4}=\{(\bv_0,\bv_2,\bv_4)\in V^3\mid \\
\bv_2\in Y_W^{ns},\ \phi(\bv_2\otimes W)\cap
\phi(\span(\bv_0,\bv_4)\otimes W)=\{\bz\}\}
\end{multline*}
and look at the map $\begin{CD}B_{0,2,4}@>{p_{0,4}}>>V^2\end{CD}$, $p_{0,4}(\bv_0,\bv_2,\bv_4)=(\bv_0,\bv_4).$

By Lemma \ref{SecondLemma}, $p_{0,4}$ is a surjection and a submersion.  Our map
$(f_0,f_4)$ actually lands in the open subset $A_{0,4}\subset V^2$, and we want to lift it
to the open subset $A_{0,2,4}\subset B_{0,2,4}$.  The restriction 
$p_{0,4}\colon A_{0,2,4}\to A_{0,4}$
is, then,  again a submersion.  For $(\bv_0,\bv_4)\in A_{0,4} $, the fiber of this restriction
over $(\bv_0,\bv_4)$ consists of all
\[\{\bv_2\in Z\mid \phi(\bv_2\otimes W)\cap\phi(\span(\bv_0,\bv_4)\otimes W)=\{\bz\}\}.\]
That is: it consists of all vectors in $V$ which are in $Y_W$ but miss $X$, which is a 
set of codimension $r(i,k)\geq r(i-1,k)$ in $V$, and also miss the set of all $\bv_2$ with
$\phi(\bv_2\otimes \bw_0)=\bc_0$, $\bc_0\in\span(\bv_0,\bv_4)\otimes W$, which,
as in the proof of Proposition \ref{FirstLemma}, has codimension greater than or equal to
\begin{multline*}
2m(i,k)-2k-k=2k^2+4ik+6k+4-3k=2k^2+4ik+3k+4
\\>k^2+2ik+i+4k-3=r(i-1,k)
\end{multline*}
for all $i\geq0$, $k>0$, because $k^2+i(2k-1)-k+7>0$.

By the induction hypothesis, since $W$ is $m(i,k)$-admissible, and therefore also
$m(i-1,k)$-admissible, on
$V$, the set of all $\bv_2\in Y_W$ which miss these two sets of codimension $\geq r(i-1,k)$
is $(i-1)$-connected.  

If we had a map of compact manifolds
$p\colon E\to B$ which was a submersion, it would be a fiber bundle.  In that case, the
condition that the fibers be $(i-1)$-connected implies, by the long exact sequence of a fibration, 
that the map $\begin{CD}\pi_i(E)@>{p_*}>>\pi_i(B)\end{CD}$
is surjective.  But we have open manifolds, instead, where fibers over different points
are not necessarily homotopy equivalent.  Nevertheless, we have seen that the fibers over different points are all $(i-1)$-connected.  By Theorem 1 in \cite{Smith}, if we have a submersion of
open manifolds where the inverse image of every point is $(i-1)$-connected, we still get the same
surjectivity: in our case $\begin{CD}\pi_i(A_{0,2,4})@>{p_{0,4\,*}}>> \pi_i(A_{0,4})\end{CD}$ is surjective.  The
result in \cite{Smith} actually requires that the fibers should all be {\em strongly} $(i-1)$-connected,
that is: every compact set in one of the fibers should be contained in an $(i-1)$-connected compact
subset of that fiber. But as explained there, strong $(i-1)$-connectedness is equivalent to
$(i-1)$-connectedness for manifolds of dimension $\geq i+2$.  Here the fiber has the same dimension
as $Y_W$ which (by Corollary \ref{DimYW}) has codimension $k$.  So
\[\dim Y_W\geq2m(i,k)-k=2k^2+4ik+5k+4>i+2,\]
and the requirement that the fibers' dimension be at least $i+2$ poses no problem for any $k> 0$, and any $(f_0,f_4)\colon S^i\to A_{0,4}$ is homotopic to
a map which can be lifted to a map $S^i\to A_{0,2,4}$.  

\smallskip

The next step is to lift a map $S^i\to A_{0,2,4}$, up to homotopy, to $A_{0,1,2,4}$.  We let
\begin{multline*}
B_{0,1,2,4}
=\{(\bv_0.\bv_1,\bv_2,\bv_4)\in V^4\mid \ 
\\
\{\bv_0,\bv_1,\bv_2\}\ {\rm are}\ W{\rm-independent},\ 
\bv_1\perp_W \bv_i,\ 0\leq i\leq 2\}
\end{multline*}
and look at the map $\begin{CD}B_{0,1,2,4}@>{p_{0,2,4}}>>V^3\end{CD}$, which omits the $\bv_1$-coordinate. Note that if we exchange the roles of $\bv_0$ and $\bv_1$,
$p_{0,2,4}$ is exactly the map $X_2\to V^2$ of (\ref{XY}) crossed with an additional copy
of $V$ (corresponding to $\bv_4$).  
Now $m(i,k)=k^2+2ik+3k+2> \frac{k^2+6k+5k-3}{2}$ since $k^2+4ik-5k+7\geq k^2-5k+7>0$, so
we have the admissibility required by Proposition \ref{FirstLemma} in the case $n=2$.
By that proposition, then, $p_{0,2,4}$ is a
submersion, and the same is true for its restriction $p_{0,2,4}\colon A_{0,1,2,4}\to A_{0,2,4}$.

Over each $(\bv_0,\bv_2,\bv_4)\in A_{0,2,4}$, the fiber $p_{0,2,4}^{-1}(\bv_0,\bv_2,\bv_4)$
consists of all
\begin{multline*}
\{\bv_1\in Z\mid \{\bv_0,\bv_1,\bv_2\}\ {\rm are}\ W{\rm -independent},\ \bv_1\perp_W
\bv_i,\ 0\leq i\leq 2,
\\
{\rm the\ lines\ from\ } \bv_1\ {\rm to}\ \bv_0\ {\rm and}\ \bv_2\ {\rm miss}\ X\}.
\end{multline*}
In other words, we have a subspace
\[V'=\{\bv\in V\mid \bv\perp_W\bv_0,\ \bv\perp_W\bv_2\}\subset V\]
of codimension $\leq 2k$ in $V$, on which $W$ is therefore by Lemma \ref{ReducingLemma}
at least $(m(i,k)-2k)$-admissible, and $m(i,k)-2k= m(i-1,k)$.  In this subspace,
we look for $\bv_1\in V'$ which is in $Y_W$ for which $\{\bv_0,\bv_1,\bv_2\}$ are $W$-independent, and which misses the cone from $\bv_0$ to $X$ and the cone from $\bv_2$ to
$X$.  The codimension of each of these cones in $V$ is at least $r(i,k)-1$, so their
codimension in $V'$ is at least $r(i,k)-2k-1=r(i-1,k)$.  And the codimension of the set of $\bv_1$ for
which $\phi(\bv_1\otimes W)\cap\phi(\span(\bv_0,\bv_2)\otimes W)\neq \{\bz\}$ is, as discussed
in the proof of Proposition \ref{FirstLemma}, at least $2m(i,k)-3k$ in $V$, so 
the codimension of the intersection of this set with $V'$ is at least $2m(i,k)-5k$ in $V'$.
We have
\[2m(i,k)-5k=2k^2+4ik+k+4>k^2+2ik+i+4k-3=r(i-1,k)\]
for $i\geq0$, $k>0$.
So by the induction hypothesis, the fibers of $p_{0,2,4}$ are $(i-1)$-connected over every point,
and as before by \cite{Smith} this means that the map
$\begin{CD}\pi_i(A_{0,1,2,4})@>{p_{0,2,4\,*}}>>\pi_i(A_{0,2,4})\end{CD}$
is surjective.

\smallskip
We can lift from $A_{0,1,2,4}$ to $A_{0,1,2,3,4}$ in the same way, using the map $X_2\to V^2$ of (\ref{XY}) crossed with two additional copies
of $V$ (corresponding to $\bv_0$ and $\bv_1$).  This finishes the lifting process and therefore
the proof.
\end{proof}


\begin{thebibliography}{KMM}
\bibitem[GM]{GM}Goresky, M. and MacPherson, R.: Stratified
Morse Theory. Springer Verlag, Berlin-Heidelberg, 1988.
\bibitem[HS]{HS}Harris, Joe; Starr, Jason: Rational curves on hypersurfaces of low degree. II.  \textit{Compos. Math.} \textbf{141}  (2005),  no. 1, 35--92.
\bibitem[H]{Har}Hartshorne, Robin: Algebraic Geometry. Graduate Texts in Mathematics, No. 52. Springer-Verlag, New York-Heidelberg, 1977. 
\bibitem[IL]{IL}Im, B. and Larsen, M.: Weak approximation for linear systems of quadrics, 
arXiv: math.NT/0306265, to appear in \textit{J. Number Theory}.
\bibitem[KMM]{KMM}Koll\'ar, J\'anos; Miyaoka, Yoichi; Mori, Shigefumi:  
Rationally connected varieties.  \textit{J. Algebraic Geom.} \textbf{1}  (1992),  no. 3, 429--448.
\bibitem[R]{Ro}Rockafellar, R. Tyrrell: Convex Analysis. Princeton Mathematical Series, No. 28 Princeton University Press, Princeton, N.J. 1970.
\bibitem[S]{Smith}Smith, J. Wolfgang: Submersions with $p$-connected fibers.  
\textit{Math. Z.} \textbf{121}  (1971), 288--294.
\end{thebibliography}
\end{document}